\documentclass[journal]{IEEEtran}

\ifCLASSINFOpdf
  \usepackage[pdftex]{graphicx}
  \graphicspath{{./Figures/}{./photos/}}
  \DeclareGraphicsExtensions{.pdf,.png}
\else
  \usepackage[dvips]{graphicx}
  \graphicspath{{./Figures/}{./photos/}}
  \DeclareGraphicsExtensions{.pdf,.png}
\fi

\usepackage{amsmath}
\ifCLASSOPTIONcompsoc
 \usepackage[caption=false,font=normalsize,labelfont=sf,textfont=sf]{subfig}
\else
 \usepackage[caption=false,font=footnotesize]{subfig}
\fi
\usepackage{dblfloatfix}
\usepackage{url}

\usepackage{amsfonts}

\usepackage{flushend}

\usepackage{listings}

\begin{document}
%
\title{Sphere Under Advection and Mean Curvature Flow}
%
%
%


\author{
  Bryce~A.~Besler,
  Tannis~D.~Kemp,
  Nils~D.~Forkert,
  Steven~K.~Boyd

\thanks{This work was supported by the Natural Sciences and Engineering Research Council (NSERC) of Canada, grant RGPIN-2019-04135.}
\thanks{B.A. Besler and T. D. Kemp are in the McCaig Institute for Bone and Joint Health, University of Calgary Canada. N.D. Forkert is with the Department of Radiology and the Hotchkiss Brain Institute, University of Calgary, Canada. S.K. Boyd is with the Department of Radiology and McCaig Institute for Bone and Joint Health, University of Calgary Canada e-mail: skboyd@ucalgary.ca}
}

%
%

\markboth{}%
{Besler \MakeLowercase{\textit{et al.}}: Sphere under Advection and Mean Curvature Flow}
%



\maketitle

\begin{abstract}
Advection and mean curvature flow is used as a model of bone microarchitecture adaptation.
It is an equivalent geometric flow to prescribed mean curvature flow with an additional rate term.
In order to validate numerical methods for simulating this flow and developing an inverse solver, a closed-form solution for advection and mean curvature flow of a sphere is derived.
\end{abstract}


\section{Introduction}
Advection and mean curvature flow has recently been used as a model for predicting bone microarchitecture changes during aging~\cite{besler2018bone}.
Towards validating numerical solver and inverse problems therein, a closed-form solution to advection and mean curvature flow of a sphere is sought.
The contents of this article are largely pedagogical.

\section{Advection and Mean Curvature Flow}
Consider an orientable, closed two-dimensional surface immersed in three dimensions $M \colon \!R^2 \rightarrow \!R^3$ with mean curvature $\kappa$ and unit normal $\hat{n}$.
A combination of advection and mean curvature flow is considered where the advection is given by a scalar rate $a$ along the unit normal direction and mean curvature is given by a rate constant $b$.
\begin{equation}
  \label{eqn:advection-mean}
  \frac{\partial M}{\partial t} = (a - b \kappa) \hat{n}
\end{equation}

Such a flow is equivalent to flow under prescribed mean curvature with an additional rate term:
\begin{equation}
  \frac{\partial M}{\partial t} = b(\tilde{\kappa} - \kappa) \hat{n}
\end{equation}
where $\tilde{\kappa} = a/b$ is the prescribed mean curvature.
The study of this flow originates from the geometric similarities between triply period minimal surfaces~\cite{schoen1970infinite,anderson1987periodic,chopp1993flow} and bone microarchitecture~\cite{hildebrand1999direct}.
It should be noted that $\tilde{\kappa}$ is not the total mean curvature and that this is not a volume preserving flow.
More precisely, the flow permits the development of singularities, producing a change in topology.

In general, $a$ can be any real number while $b$ should be a real number greater than or equal to zero.
A negative value for $b$ would imply inverse mean curvature flow, which is unstable when points on the surface have zero mean curvature.
One can see that the flow stops when $\kappa = \tilde{\kappa} = a/b$ everywhere on the surface.

\section{The Sphere}
This paper is concerned with a closed-form solution to a sphere moving under advection and mean curvature.
All work will be done on the two-sphere $S_2$.

Consider the two-sphere mapping to spherical coordinates $M \colon (u,v) \rightarrow (\rho, \theta, \phi)$ with radius $r_0$.
\begin{align}
  \label{eqn:sphere:mapping}
  \begin{pmatrix}
    \rho \\
    \theta \\
    \phi
  \end{pmatrix} = 
  \begin{pmatrix}
    r_0 \cos u \sin v \\
    r_0 \sin u \sin v \\
    r_0 \cos v
  \end{pmatrix}
\end{align}
The normal is along the radial direction.
\begin{equation}
  \label{eqn:sphere:normal}
  \hat{n}(u,v) = 1\hat{\rho} + 0 \hat{u} + 0 \hat{v}
\end{equation}
Since the surface normal is aligned with the $\hat{\rho}$ direction, all curve evolution problems will reduce to a differential equation on the radius by spherical symmetry.
The mean curvature at every point is the inverse of radius.
\begin{equation}
  \label{eqn:sphere:curvature}
  \kappa(u,v) = 1/r
\end{equation}
These prerequisites allow the development of a closed-form solution for the sphere.

\section{Sphere Under Advection}
\label{sec:advection}
Lets begin with the problem of advection ($b = 0$).
In this case, Equation~\ref{eqn:advection-mean} reduces to a simple expression:
\begin{equation}
  \frac{\partial M}{\partial t} = a \hat{n}
\end{equation}
Substituting Equations~\ref{eqn:sphere:mapping} and \ref{eqn:sphere:normal}:
\begin{align}
  \begin{pmatrix}
    \rho_t \\
    \theta_t \\
    \phi_t
  \end{pmatrix} = a \cdot
  \begin{pmatrix}
    1 \cdot \hat{\rho} \\
    0 \cdot \hat{\theta} \\
    0 \cdot \hat{\phi}
  \end{pmatrix}
\end{align}
where $\left(\cdot\right)_t$ is shorthand for a temporal derivative.
This gives a single initial value problem to solve:
\begin{equation}
  \left\{
    \begin{array}{ll}
      \rho_t = a\\
      \rho(0) = r_0
    \end{array}
  \right.
\end{equation}
The solution is immediate:
\begin{equation}
  \label{eqn:advection:solution}
  \rho(t) = r_0 + at
\end{equation}

\begin{figure}[t]
  \centering
    \includegraphics[width=0.9\linewidth]{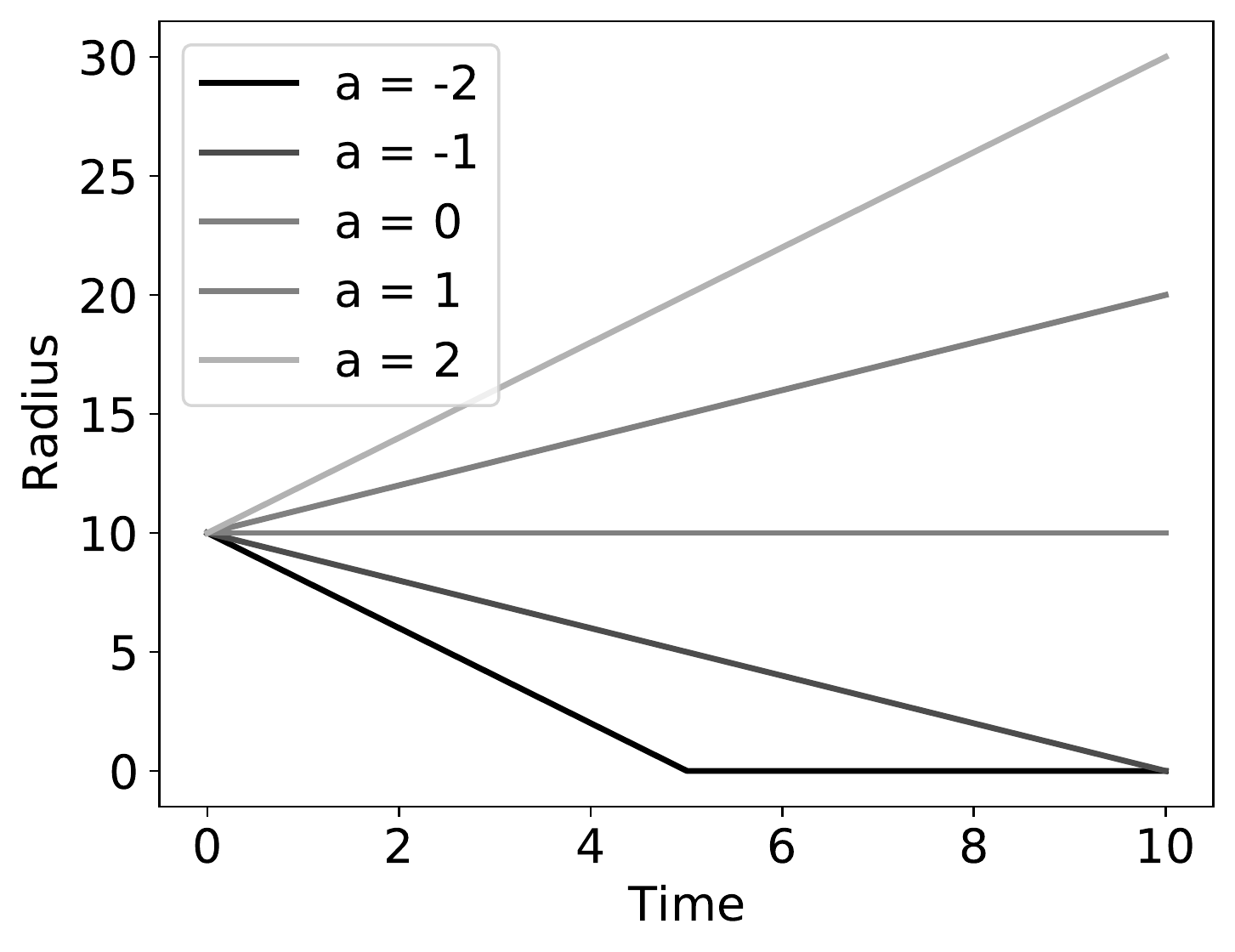}%
  \caption{Advection solution to the sphere ($r_0 = 10$).}
  \label{fig:advection}
\end{figure}

Equation~\ref{eqn:advection:solution} agrees with intuition.
The surface is moving at a linear rate of $a$ units of distance per units of time along the normal of the sphere.
If $a$ is positive, the sphere grows forever.
If $a$ is negative, it shrinks until it vanishes at the time $t = -r_0/a$.
The solution is plotted for five values of $a$ in Figure~\ref{fig:advection}.

\section{Sphere Under Mean Curvature Flow}
\label{sec:meanflow}
Attention is now placed on mean curvature flow ($a = 0$).
In this case, Equation~\ref{eqn:advection-mean} reduces to a simple expression:
\begin{equation}
  \label{eqn:mean}
  \frac{\partial M}{\partial t} = -b \kappa \hat{n}
\end{equation}
Substituting Equations~\ref{eqn:sphere:mapping}, \ref{eqn:sphere:normal} and \ref{eqn:sphere:curvature}:
\begin{align}
  \begin{pmatrix}
    \rho_t \\
    \theta_t \\
    \phi_t
  \end{pmatrix} = -b / \rho \cdot
  \begin{pmatrix}
    1 \cdot \hat{\rho} \\
    0 \cdot \hat{\theta} \\
    0 \cdot \hat{\phi}
  \end{pmatrix}
\end{align}
Again, this leads to a single initial value problem:
\begin{equation}
  \left\{
    \begin{array}{ll}
      \rho_t = -b/\rho\\
      \rho(0) = r_0
    \end{array}
  \right.
\end{equation}
Through some substitution, the problem can be solved:
\begin{equation}
  \label{eqn:mean:solution}
  \rho(t) = \sqrt{r_0^2 - 2bt}
\end{equation}
This solution is a classic pedagogical result in mean curvature flow~\cite{bellettini2014lecture}.

\begin{figure}[t]
  \centering
    \includegraphics[width=0.9\linewidth]{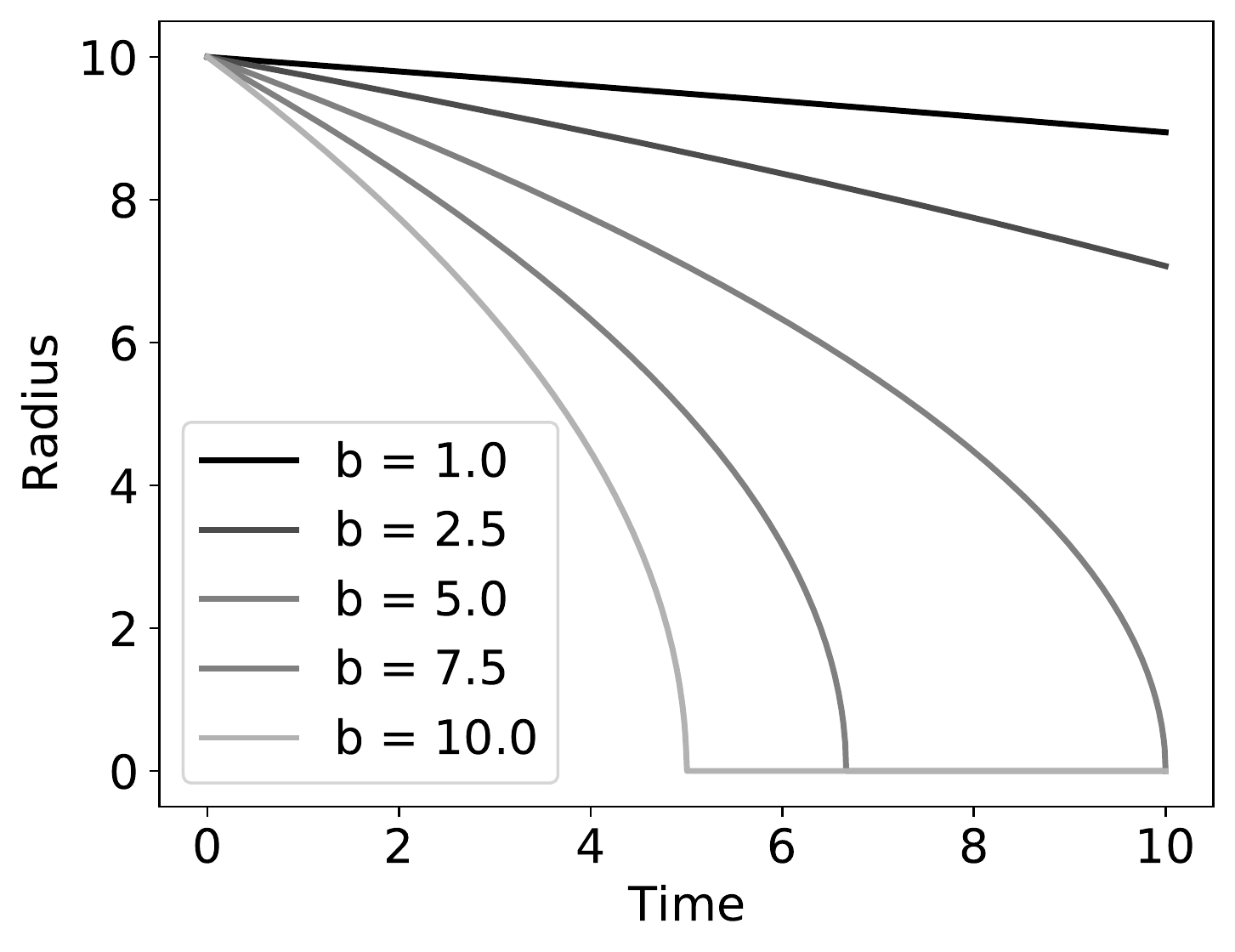}%
  \caption{Mean curvature solution to the sphere ($r_0 = 10$).}
  \label{fig:mean}
\end{figure}

As with advection, Equation~\ref{eqn:mean:solution} agrees with intuition.
Since the mean curvature everywhere on a sphere is positive and $b$ is positive, the negative sign in Equation~\ref{eqn:mean} suggests that the sphere is always shrinking.
Indeed, the sphere shrinks until it vanishes at $t = r_0^2/2b$.
The solution is plotted for five values of $b$ in Figure~\ref{fig:mean}.

\section{Sphere Under Advection and Mean Curvature Flow}
Having some background and intuition, the main solution is sought.
Skipping the middle steps as described in Sections~\ref{sec:advection} and \ref{sec:meanflow}, the initial value problem can be stated.
\begin{equation}
  \label{eqn:base}
  \left\{
    \begin{array}{ll}
      r_t = a - \frac{b}{r}\\
      r(0) = r_0
    \end{array}
  \right.
\end{equation}
Notation is changed from $\rho$ to $r$ for ease of interpretation.

\begin{figure}[b]
  \centering
    \includegraphics[width=0.9\linewidth]{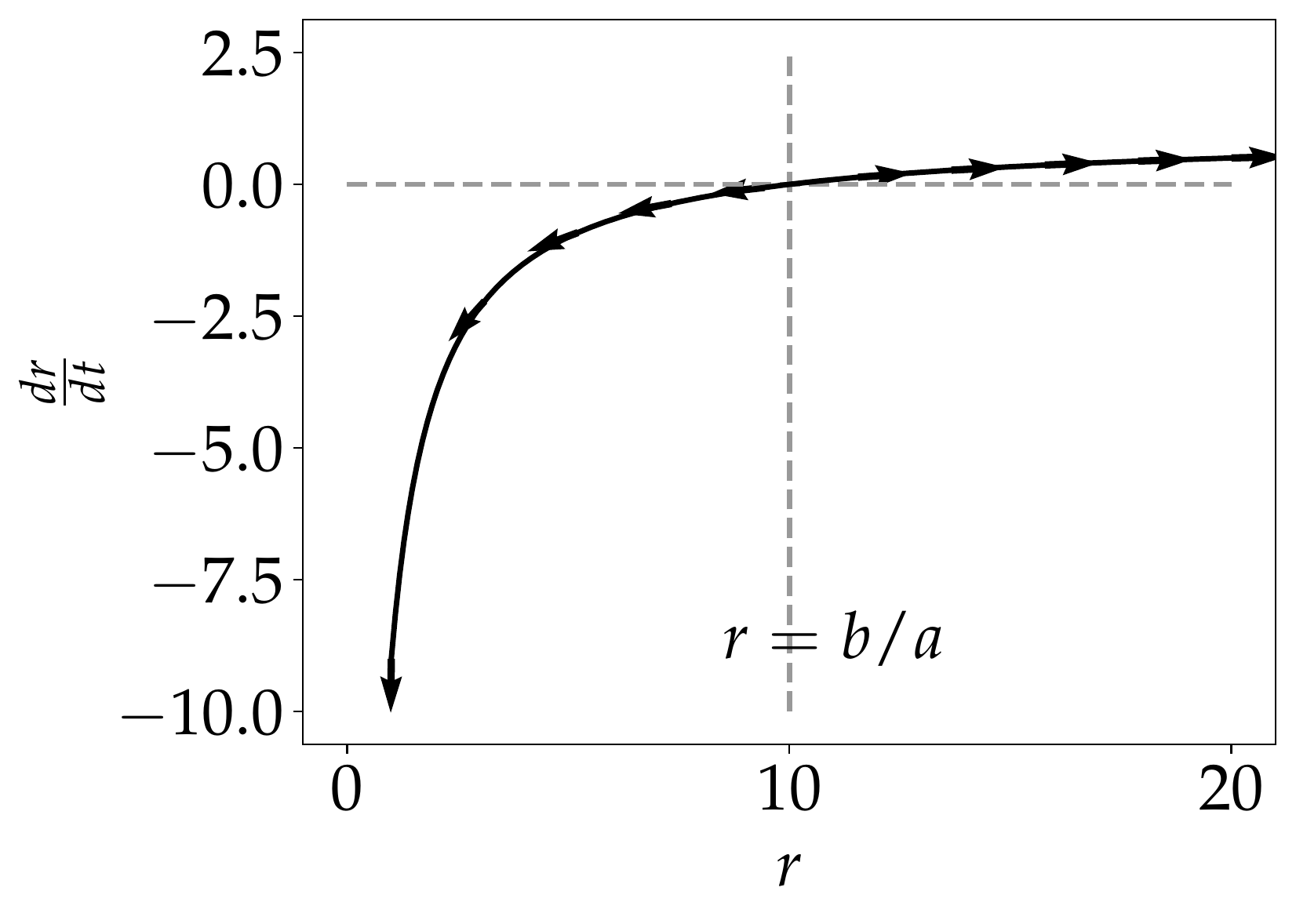}%
  \caption{Plotting the meta-stable point of the sphere under advection and mean curvature ($a = 1$, $b = 10$).}
  \label{fig:loss}
\end{figure}

Analyzing Equation~\ref{eqn:base} can give insight into the model.
In general, the model stops flowing when $r_t = 0$, which corresponds to $r = b/a$.
However, this is only the case when $a$ is positive.
When $a$ is negative, the sphere shrinks forever.
When $a$ is positive, there exists a point where the shrinking under mean curvature is balanced by the growth from advection.
However, this point is only meta-stable.
If the initial sphere radius is exactly $r_0 = b/a$, the sphere will be constant over time.
However, if the radius is slightly increased, the mean curvature term decreases and the sphere grows.
If the radius is slightly shrunk, the mean curvature term increases and the sphere shrinks.
In summary:
\begin{itemize}
  \item $a < 0 \rightarrow$ shrink to zero
  \item $0 < a < \frac{b}{r_0} \rightarrow $ shrink to zero
  \item $ a = \frac{b}{r_0} \rightarrow $ meta-stable
  \item $\frac{b}{r_0} < a \rightarrow $ grow to infinity
\end{itemize}
In general, away from the meta-stable point, when growing, growth is like advection. When shrinking, shrinking is like mean curvature flow.
This is demonstrated in a phase diagram in Figure~\ref{fig:loss}.

\begin{figure}[b]
  \centering
    \includegraphics[width=0.9\linewidth]{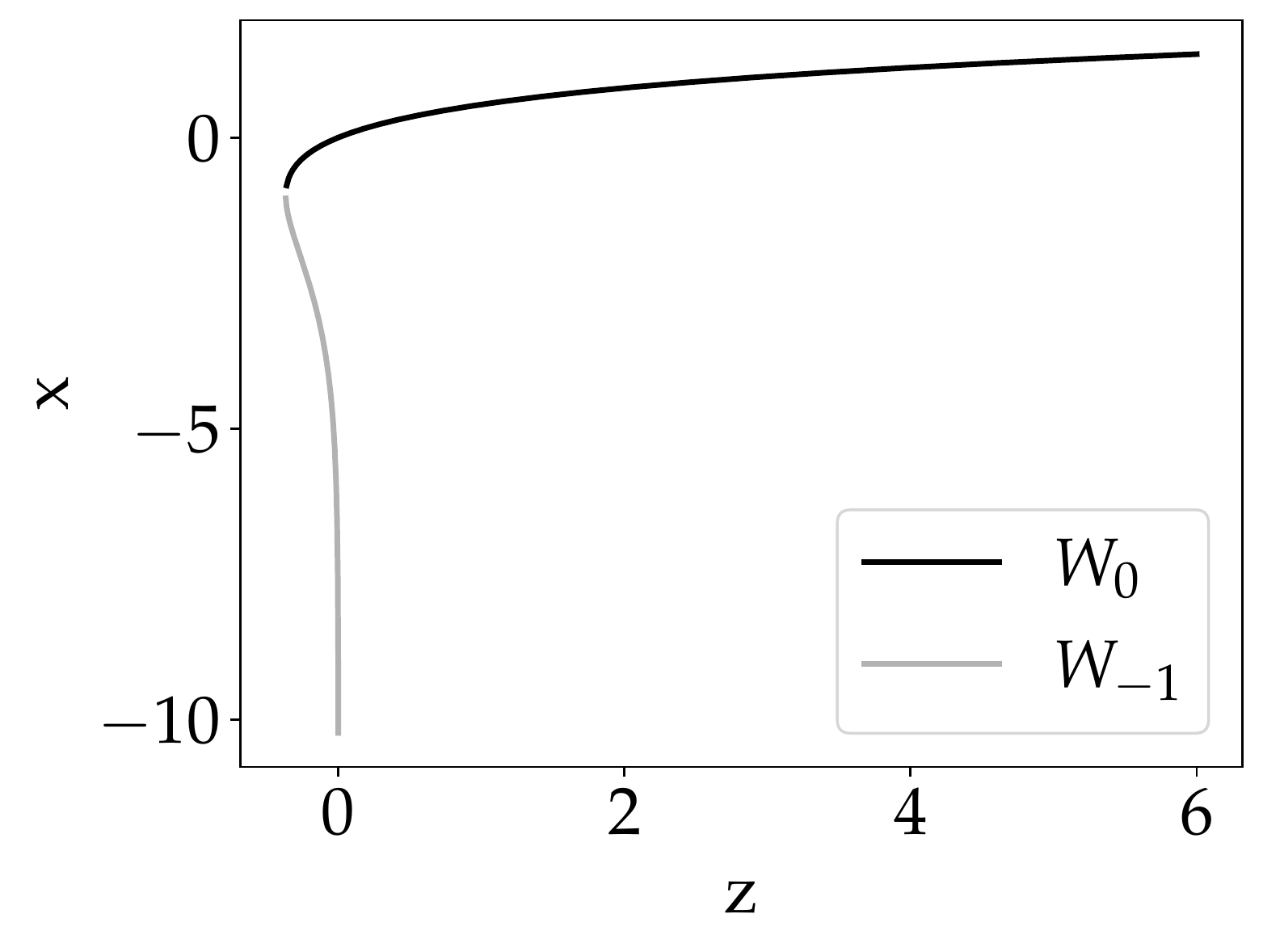}%
  \caption{Plotting the Lambert W function for real numbers.}
  \label{fig:w}
\end{figure}

Now, the closed-form solution to Equation~\ref{eqn:base} is derived.
\begin{eqnarray}
  \frac{dr}{dt} = a - \frac{b}{r} \\
  \frac{r}{ar - b} dr = dt
\end{eqnarray}
Using the substitution $x = ar-b$ leads to:
\begin{eqnarray}
  \frac{x+b}{a}\frac{1}{x} \frac{dx}{a} = dt \\
  \left[\frac{1}{a^2} + \frac{b}{ax}\right] dx = dt
\end{eqnarray}
Which can be integrated:
\begin{equation}
  \frac{ar-b}{a^2} + \frac{b}{a^2}\ln\left(ar-b\right) = t + c
\end{equation}
It is convenient to solve for $c$ at this point with the initial condition $r(0) = r_0$.
\begin{equation}
  c = \frac{ar_0-b}{a^2} + \frac{b}{a^2}\ln\left(ar_0 - b\right)
\end{equation}
The final step is to solve for $r$.
\begin{eqnarray}
  \label{eqn:near-solved}
  \exp\left( \frac{ar-b}{b}\right) \frac{ar-b}{b} = \frac{1}{b}\exp\left(\frac{a^2}{b}(t+c)\right)
\end{eqnarray}
While this appears unsolvable, the left hand side of Equation~\ref{eqn:near-solved} is known as the ``product log" or Lambert W function~\cite{lambert1758observationes,weisstein2002lambert}.
\begin{eqnarray}
  xe^x &=& z \\
  x &=& W_k(z)
\end{eqnarray}
The Lambert W function is plotted in Figure~\ref{fig:w}.
In the case of Equation~\ref{eqn:near-solved},
\begin{eqnarray}
  \label{eqn:x}
  x &=& \frac{ar-b}{b} \\
  z &=& \frac{1}{b} \exp\left(\frac{a^2}{b}(t+c)\right)
\end{eqnarray}

Using the Lambert W and substituting $c$ allows us to solve for $r$:
\begin{eqnarray}
  \frac{ar-b}{b} = W_k\left(\frac{1}{b}\exp\left(\frac{a^2}{b}(t+c)\right)\right) \\
  \label{eqn:soln}
  r = \frac{b}{a}\left(W_k\left[\frac{ar_0-b}{b} \exp\left(\frac{ar_0-b}{b}\right)\exp\left(\frac{a^2t}{b}\right)\right] + 1\right)
\end{eqnarray}

Next, the appropriate branch of $W_k$ must be selected since multiple solutions are possible.
Since this problem deals with real and not complex numbers, only $W_0$ and $W_{-1}$ are available.
From Figure~\ref{fig:w}, two solutions can be seen for $z < 0$, the solutions changing when $x = -1$.
Using Equation~\ref{eqn:x}, the switching condition can be defined.
\begin{eqnarray}
  \frac{ar-b}{b} < -1 \\
  ar < 0
\end{eqnarray}
However, $r$ is always positive since it is the radius of a sphere.
As such, use $W_{-1}$ if $a<0$ and use $W_0$ if $a>0$.
Note that this piecewise function is continuous at $a=0$.

The final step is to derive a vanishing time for the solution.
Again, the vanishing time is the $t$ when $r(t)=0$, which is when $W_k(z) = -1$.
This corresponds to the situation when the argument of $W_k$ in Equation~\ref{eqn:soln} is equal to $\frac{-1}{e}$.
Looking at that argument, the vanishing time is found.
\begin{eqnarray}
  \frac{ar_0-b}{b} \exp\left(\frac{ar_0-b}{b}\right)\exp\left(\frac{a^2t}{b}\right) = \frac{-1}{e} \\
  \label{eqn:vanish1}
  t = \frac{b}{a^2} \ln\left(\frac{b}{b-ar_0}\right) - \frac{r_0}{a}
\end{eqnarray}
Equation~\ref{eqn:vanish1} will only be valid when the argument to the natural logarithm is positive.
\begin{eqnarray}
  \frac{b}{b-ar_0} > 0 \\
  b - ar_0 > 0 \\
  r_0 < \frac{b}{a}
\end{eqnarray}
This is the same condition that was found earlier from analysis of the differential equation.
The equation for vanishing time is then given:
\begin{equation}
  \label{eqn:vanish}
  t = \left\{
    \begin{matrix}
      \infty & \text{ if } r_0 \geq \frac{b}{a} \\
      \frac{b}{a^2} \ln\left(\frac{b}{b-ar_0}\right) - \frac{r_0}{a} & \text{ if } r_0 < \frac{b}{a}
    \end{matrix}
  \right.
\end{equation}
Taking the limits as $a \rightarrow 0$ or as $b \rightarrow 0$ of Equation~\ref{eqn:vanish}, the advection and mean curvature vanishing times from before can be found.

\begin{figure}[t]
  \centering
  \begin{tabular}{c}
    \subfloat[$r_0 = 11, a = 1, b = 10$]{
      \includegraphics[width=0.8\linewidth]{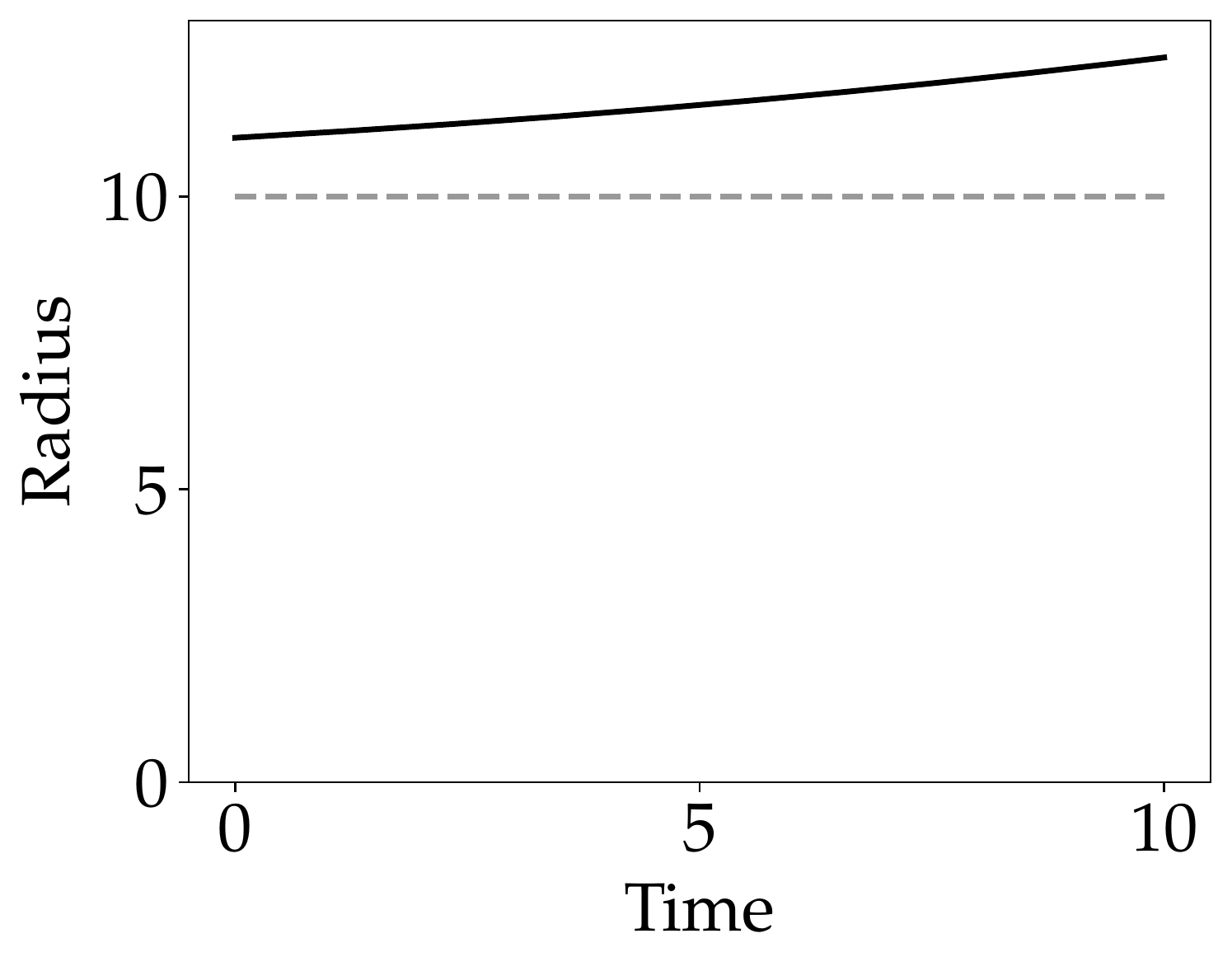}%
      \label{fig:solution:1}
    } \\
    \subfloat[$r_0 = 9, a = 1, b = 10$]{
      \includegraphics[width=0.8\linewidth]{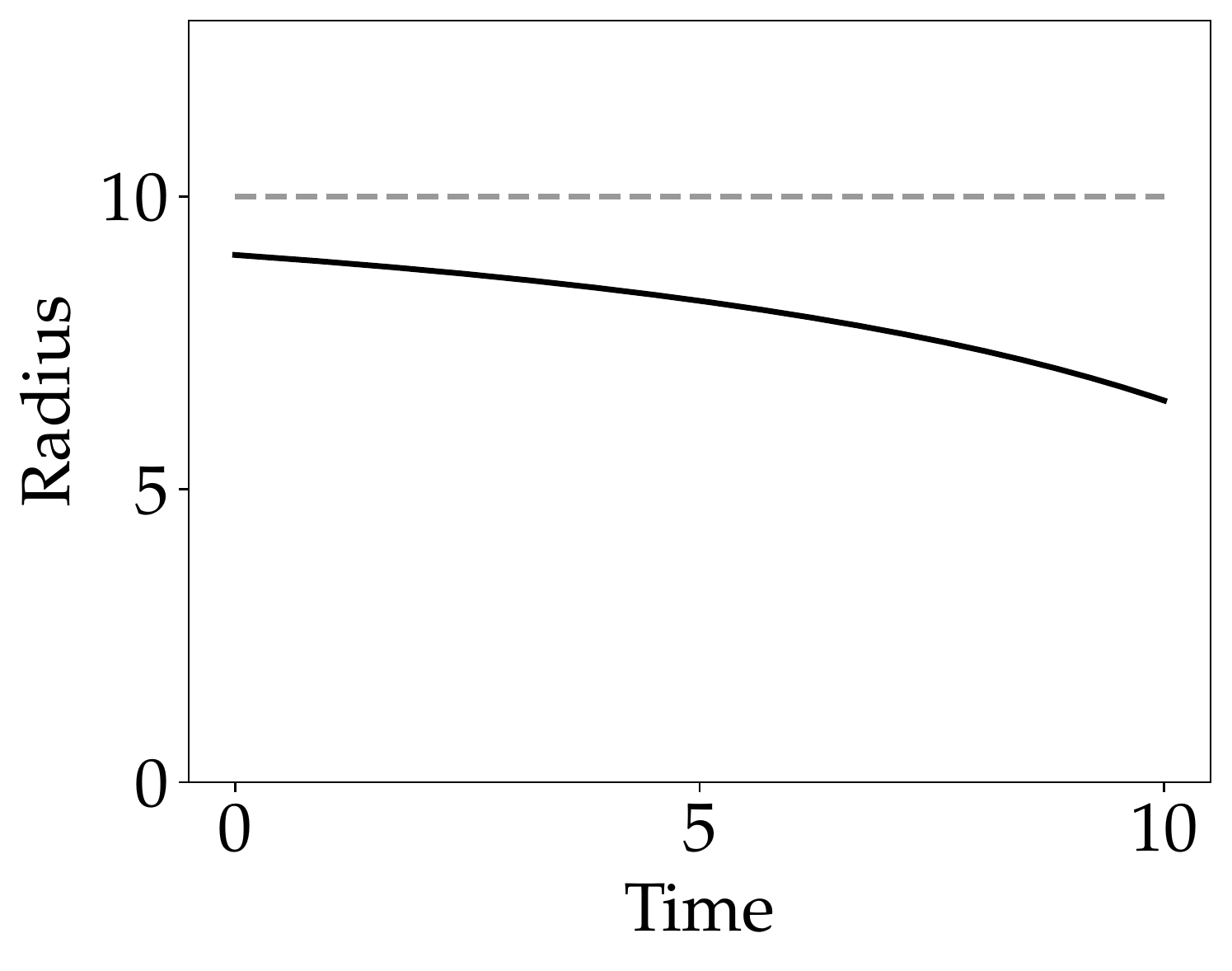}%
      \label{fig:solution:2}
    } \\
    \subfloat[$r_0 = 10, a = -1, b = 10$]{
      \includegraphics[width=0.8\linewidth]{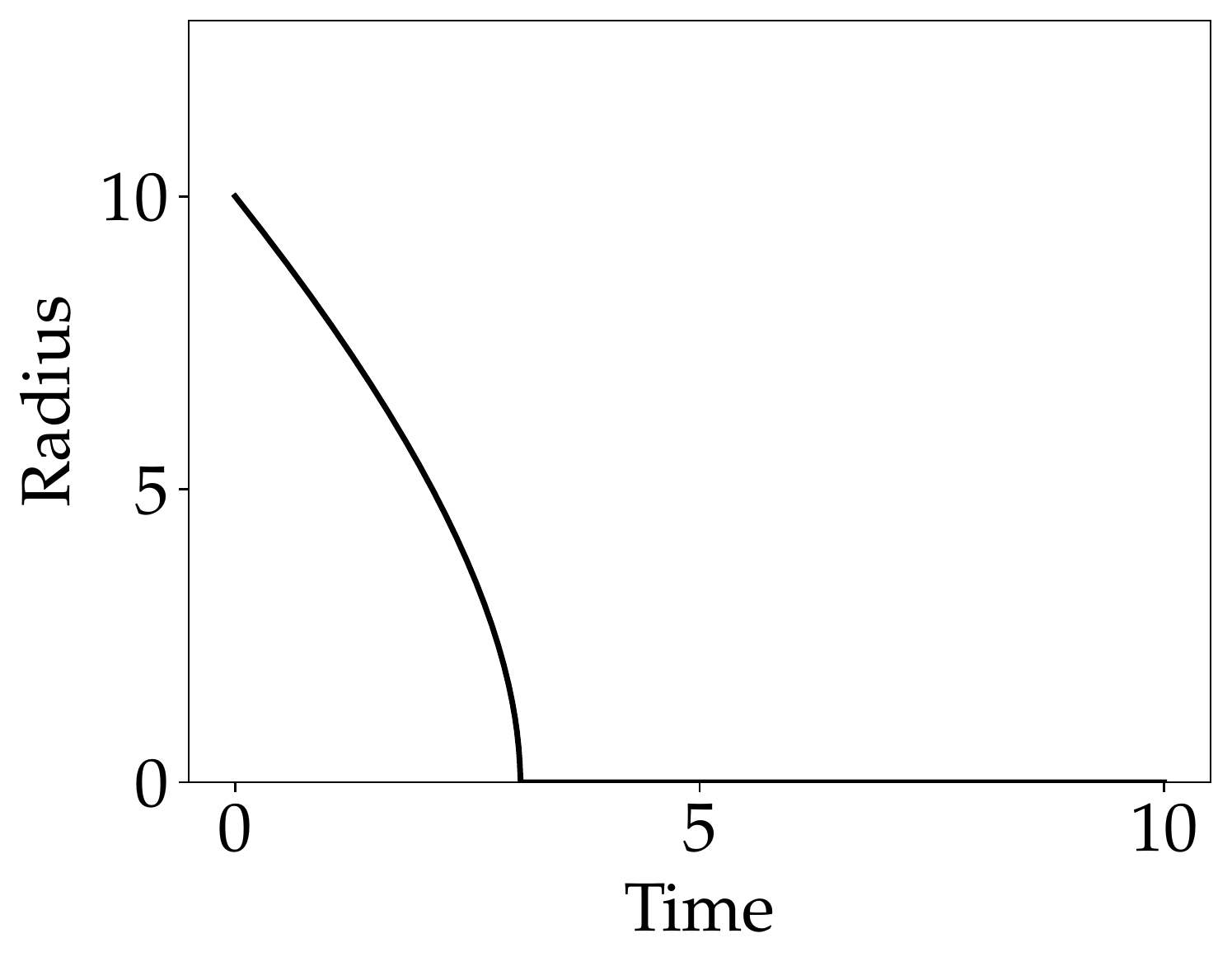}%
      \label{fig:solution:3}
    }
  \end{tabular}
  \caption{Solutions to the advection and mean curvature flow of a sphere. The three conditions correspond to (\ref{fig:solution:1}) growth, (\ref{fig:solution:2}) shrinking, and (\ref{fig:solution:1}) shrinking with negative advection. The dashed line denotes the meta-stable point.}
  \label{fig:solution}
\end{figure}

An implementation of the closed-form solution is given in the Appendix.
To finalize the analysis, the solution is plotted in Figure~\ref{fig:solution} for various parameters.
The solution is only of interest around the meta-stable point.
Outside the meta-stable point, one of the two terms is much smaller than the other, meaning only one term drives the flow.
This is important to consider for inverse solvers which may not be able to accurately solve for the non-dominant parameter in this simple geometry.

\section{Conclusion}
A closed-form solution for the motion of a sphere under advection and mean curvature flow is developed.
Vanishing times for the sphere are also calculated.
This model provides a closed-form solution for testing numerical solvers of advection and mean curvature flow and inverse problems therein.


%


\appendices

\appendix[Source Code]
\label{app:python}
Source code is provided in Python for ease of implementation.

\lstdefinestyle{CustomStyle}{
  language=Python,
  numbers=none,
  stepnumber=1,
  numbersep=10pt,
  tabsize=4,
  showspaces=false,
  showstringspaces=false,
  breaklines=true
}
\lstset{basicstyle=\tiny,style=CustomStyle}
\begin{lstlisting}
# Imports
import numpy as np
from scipy.special import lambertw

def sphere_under_advection(t, r_0, a):
    '''Compute the radius of a sphere under advection

    Parameters
    ----------
    t : float, np.array
        Time(s) for which to solve
    r_0 : float
        Initial radius of the sphere
    a : float
        Advection constant (positive = grow)
    '''
    # Compute r
    r = a*t + r_0
    
    # Check if disappear
    if np.abs(a) > np.finfo(float).eps:
        t_vanish = -r_0 / a
        if t_vanish > 0:
            r[t >= t_vanish] = 0.0
    
    return r

def sphere_under_mean_curvature(t, r_0, b):
    '''Compute the radius of a sphere under advection

    Parameters
    ----------
    t : float, np.array
        Time(s) for which to solve
    r_0 : float
        Initial radius of the sphere
    b : float
        Mean curvature constant (must be > 0)
    '''
    # Compute vanishing time
    t_vanish = np.Inf
    if np.abs(b) > np.finfo(float).eps:
        t_vanish = r_0**2/(2.0*b)

    # Split t by t_vanish to avoid sqrt(-1)
    t_pos = t[t<t_vanish]
    t_neg = t[t>=t_vanish]

    # Compute r
    r = np.concatenate([
        np.sqrt(r_0**2 - 2*b*t_pos),
        np.zeros_like(t_neg)
    ])

    return r

def sphere_vanish_time(r_0, a, b):
    '''Compute the vanishing time of a sphere under advection and mean curvature flow

    Parameters
    ----------
    r_0 : float
        Initial radius of the sphere
    a : float
        Advection constant (positive = grow)
    b : float
        Mean curvature constant (must be > 0)
    '''
    # Switch on infinity
    cond = a*r_0 < b
    
    t_vanish = np.Inf
    if cond:
        t_vanish = b/a**2 * np.log(b / (b - a *r_0)) - r_0/a
    
    return t_vanish

def sphere_under_advection_and_mean_curvature(t, r_0, a, b):
    '''Compute the vanishing time of a sphere under advection and mean curvature flow

    Parameters
    ----------
    t : float, np.array
        Time(s) for which to solve
    r_0 : float
        Initial radius of the sphere
    a : float
        Advection constant (positive = grow)
    b : float
        Mean curvature constant (must be > 0)
    '''

    # Run simpler methods if possible
    if b == 0:
        return sphere_under_advection(t, r_0, a)
    if a == 0:
        return sphere_under_mean_curvature(t, r_0, b)
    
    # Determine vanishing time
    t_vanish = sphere_vanish_time(r_0, a, b)
    
    # Select k branch of W_k
    k = -1 if a<0 else 0
    
    # Compute our values
    t_not_vanish = t[t<t_vanish]
    x = (a*r_0-b)/b
    z = x * np.exp(x) * np.exp(a**2*t_not_vanish/b)
    w_ = np.real(lambertw(z, k=k))
    r_not_vanish = (b/a) * (w_ + 1.0)
    
    # Combine with vanish
    r = np.concatenate([
        r_not_vanish,
        np.zeros_like(t[t>=t_vanish])
    ])
    
    return r
\end{lstlisting}



\ifCLASSOPTIONcaptionsoff
  \newpage
\fi



\bibliographystyle{IEEEtran}
\bibliography{IEEEabrv,sphere.bib}

\end{document}